\documentclass[11pt]{amsart}
\usepackage{amsmath,amssymb,amsfonts}

\author{Christine Bachoc}
\address{C. Bachoc, Laboratoire A2X, Universit\'e Bordeaux I,
 351, cours de la Li\-b\'e\-ra\-tion, 33405 Talence
France }
\email{bachoc@math.u-bordeaux1.fr}


\title{Designs, groups and lattices}

\date{\today} 

\newtheorem{theorem}{Theorem}[section]

\newtheorem{remark}[theorem]{Remark}
\newtheorem{corollary}[theorem]{Corollary}
\newtheorem{lemma}[theorem]{Lemma}
\newtheorem{definition}[theorem]{Definition}

\newcommand{\Z}{{\mathbb{Z}}}

\newcommand{\R}{{\mathbb{R}}}
\newcommand{\C}{{\mathbb{C}}}
\newcommand{\F}{{\mathbb{F}}}

\newcommand{\Ck}{{\mathcal{C}}_k}
\newcommand{\Gk}{{\mathcal{G}}_k}
\newcommand{\Pk}{{\mathcal{P}}_k}
\newcommand{\G}{{\mathcal{G}}}

\newcommand{\Gmn}{{\mathcal{G}_{m,n}}}

\newcommand{\On}{{\operatorname{O}(\R^n)}}

\newcommand{\Hmu}{\operatorname{H}_{m,n}^{2\mu}}

\newcommand{\D}{{\mathcal{D}}}
\newcommand{\Dsig}{{\mathcal{D}}_{\Sigma}}

\newcommand{\Sm}{\operatorname{S}_m(\R)}

\newcommand{\Smsd}{\operatorname{S}_m^{\geq 0}(\R)}
\newcommand{\Smd}{\operatorname{S}_m^{>0}(\R)}

\newcommand{\Aut}{\operatorname{Aut}}

\newcommand{\End}{\operatorname{End}}

\newcommand{\pr}{\operatorname{pr}}

\newcommand{\eproof}{\hfill $\fbox{}$ \\}

\begin{document}

\maketitle    

\begin{abstract}
The notion of designs in Grassmannian spaces was introduced by the author
and R. Coulangeon, G. Nebe, in \cite{BCN}. After having recalled 
some basic properties of these objects and the connections with the theory
of lattices, we prove that the sequence of Barnes-Wall lattices hold
$6$-Grassmannian designs. We also discuss the connections between 
the notion of Grassmannian design and the notion of design associated with
the symmetric space of the totally isotropic subspaces 
in a binary quadratic space, which is revealed in a certain construction
involving the Clifford group.
\end{abstract}

\section{Introduction}

Roughly speaking, a design is a finite subset of a space $X$ which
``approximates well'' $X$. In the case of finite spaces $X$, such
objects arose from different contexts like statistics, finite geometries,
graphs, and are well understood in the framework of association schemes
(\cite{BI}, \cite{DL}).
Later the notion of designs was extended to the two-point homogeneous
real manifolds (\cite{DGS}). Of special interest are the so-called spherical
designs, defined on the unit sphere of the Euclidean space. 
Due mainly to the work of Boris Venkov (\cite{V1}), we know that nice spherical designs 
arise from certain families of lattices, and that the lattices which contain
spherical designs are locally dense. Moreover,
this combinatorial property gives a hint to classify these lattices, which
was recently fulfilled in many cases (\cite{BV}).

In a common work with R. Coulangeon and G. Nebe, we have generalized
these notions to the real Grassmannian  spaces $\Gmn$. 
This was the subject of my talk at the 
XXIII\`emes Journ\'ees Arithm\'etiques (2003), in Graz.
I have chosen not to reproduce this talk here, but rather to 
present some complementary results on one aspect of this subject,
which was not emphasized in Graz, namely the links with group representation.
In particular, we will not discuss at all here the connections with
Siegel modular forms.

Sections 2 to 4 essentially review on results from \cite{BCN}.
The zonal polynomials associated
with the action of the orthogonal group on $\Gmn$, which are generalized Jacobi
polynomials in $m$ variables, play a crucial role. They are presented in \S 2.
The existence of Grassmannian designs in a lattice  is connected to its
Rankin functions $\gamma_{m,n}$, this is recalled in \S 3. 
In \S 4, we recall how certain  properties of the representations of a 
finite subgroup of $\On$ ensures that its orbits on $\Gmn$ are designs.
This is successfully applied to the automorphism group of many lattices.

In \S 5, we introduce the Clifford groups $\Ck<O(\R^{2^k})$,
and their subgroups $\G_k$, of index $2$,  which are the automorphism groups
of the Barnes-Wall lattices. These groups have recently attracted
attention in combinatorics because of their appearance 
in several apparently disconnected situations (finite geometries,
quantum codes, lattices, Kerdock codes..). In \cite{NRS}, a very nice
combinatorial proof that their polynomial invariants are spanned by the
generalized weight enumerators of binary codes is given.
We partly extend this result to the subgroup $\G_k$.
As a consequence, we obtain that 
the Barnes-Wall lattices support
Grassmannian $6$-designs, and that they are local maxima for all the 
Rankin constants.

The last section, \S 6, discusses some other  constructions of 
Grassmannian designs associated with the Clifford groups. 
We encounter another notion of design, this time associated
with the space of totally isotropic subspaces of fixed dimension
in a binary quadratic space. This space is homogeneous and 
symmetric for the action of the corresponding binary orthogonal group. 
Unsurprisingly, the Clifford group connects these two notions of designs,
leading to interesting new examples of Grassmannian designs.

\section{Grassmannian designs}

\subsection{Definitions}

The notion of Grassmannian designs was introduced in 
\cite{BCN}. 
Let $m\leq n/2$, and let $\Gmn$ be the real Grassmannian space, together with the transitive
action of the real orthogonal group $\On$. The starting point is the
decomposition of the space of complex-valued squared module integrable functions
$L^2(\Gmn)$ under the action of $\On$. One has:

\begin{equation}\label{e1}
L^2(\Gmn)=\oplus_{\mu}\Hmu
\end{equation}
where the sum is over the partitions $\mu=\mu_1\geq \dots\geq \mu_m\geq 0$, and the spaces $\Hmu$ are isomorphic to
the irreducible representation of $\On$ canonically associated with
$2\mu$, and denoted $V_n^{2\mu}$ (see \cite{GW}).
Here  $2\mu=2\mu_1\geq \dots\geq 2\mu_m\geq 0$ is a partition with even parts.
The degree of the partition 
$\mu$ is by definition $\deg(\mu):=\sum_{i} \mu_i$ and its length
$l(\mu)$ is the number of its non-zero parts.

As an example, when $l(\mu)=1$, the representation $V_n^{\mu}$ is isomorphic
to the space of polynomials in $n$ variables, homogeneous of degree 
$\mu_1$, and harmonic, i.e. annihilated by the standard Laplace operator.
When $l(\mu)>1$, the representations 
$V_n^{\mu}$ have more complicated but still explicit realizations
as spaces of polynomials in matrix arguments.

\begin{definition}[\cite{BCN}]\label{def1}
A finite subset $\mathcal D$ of $\Gmn$ is called a $2t$-design if,
for all $f\in \Hmu$ and all $\mu$ with $0\leq \deg(\mu)\leq t$,
\begin{equation}
\int_{\Gmn}f(p)dp=\frac{1}{|\mathcal D|}\sum_{x\in \mathcal D}f(x).
\end{equation}
\end{definition}

The decomposition (\ref{e1}) immediately shows that this definition is 
equivalent to the condition:

\begin{equation}\label{e2}
\text{for all } f\in \Hmu \text{ and all }\mu\text{ with }1\leq \deg(\mu)\leq t,
\sum_{x\in \mathcal D}f(x)=0.
\end{equation}

There is a nice characterization of the designs in terms of the 
zonal functions of $\Gmn$, which is much more satisfactory from
the algorithmic point of view. We briefly recall it here.

It is a classical fact that the 
orbits under the action of $\On$ of the pairs $(p,p')$ of
elements of $\Gmn$ are characterized by their  so-called principal angles
$(\theta_1,\dots,\theta_m)\in [0,\pi/2]^m$. We set
$y_i:=\cos^2(\theta_i)$. The polynomial functions on
$\Gmn\times\Gmn$ which are invariant under the simultaneous
action of $\On$ are polynomials in the variables $(y_1,\dots,y_m)$,
and their space is isomorphic to the algebra 
$\C[y_1,\dots,y_m]^{S_m}$ of symmetric polynomials in $m$ variables.
Moreover, there is a unique sequence of orthogonal polynomials
$P_{\mu}(y_1,\dots,y_m)$ indexed by the partitions of length $m$,
such that 
$\C[y_1,\dots,y_m]^{S_m}=\oplus_{\mu}\C P_{\mu}$,
$P_{\mu}(1,\dots,1)=1$, and the function~:
$p\in \Gmn \to P_{\mu}(y_1(p,p'), \dots, y_m(p,p'))$ defines, for all
$p' \in \Gmn$, an element of $\Hmu$. These polynomials have degree
$\deg(\mu)$. They are explicitly calculated in \cite{JC}, where it is shown that
they belong to the family of Jacobi polynomials.

More precisely, James and Constantine  show that the canonical measure 
on $\Gmn$,
induces on $\C[y_1,\dots,y_m]^{S_m}$
the following measure:

$$d\mu(y_1,\dots,y_m)=
\lambda \prod
_{1\leq i<j\leq m}|y_i-y_j|\prod_{1\leq i\leq m}y_i^{-1/2}(1-y_i)^{n/2-m-1/2} dy_i
$$
(where $\lambda$ is chosen so that $\int_{[0,1]^m}d\mu(y_1,\dots,y_m)=1$).
This measure defines an hermitian 
product on $\C[y_1\dots,y_m]^{S_m}$, namely

\begin{equation*}
[f,g]=\int_{[0,1]^m}f(y)\overline{g(y)}d\mu(y).
\end{equation*}

Since the irreducible subspaces $\Hmu$ are pairwise orthogonal, the corresponding polynomials $P_{\mu}$ must be orthogonal for this hermitian product.
Together with some knowledge on the monomials of degree $\deg(\mu)$
that occur in $P_{\mu}$, it is enough to uniquely determine them.
However, the most efficient way 
to calculate them is to exploit the fact that they are eigenvectors
for the operator on $\C[y_1,\dots,y_m]^{S_m}$ induced by the 
Laplace-Beltrami operator (see \cite{JC}, \cite{BCN} for more details). 

The first ones are equal to:

\newpage

\begin{align*}
P_0 &=1 \\
P_{(1)} &=\frac{1}{\beta_1}\left(\sum y_i -\frac{m^2}{n}\right), \quad
\beta_1=m(1-\frac{m}{n}) \\
P_{(11)} &=\frac{1}{\beta_{11}}\left(\sum y_iy_j -\frac{(m-1)^2}{n-2}\sum y_i+\frac{m^2(m-1)^2}{2(n-1)(n-2)}\right), \\
&\beta_{11}=\frac{m(m-1)}{2}(1-2\frac{m-1}{n-2}+\frac{m(m-1)}{(n-1)(n-2)})\\
P_{(2)}&=\frac{1}{\beta_2}\left(\sum y_i^2 +\frac{2}{3}\sum y_iy_j -
\frac{2(m+2)^2}{3(n+4)}\sum y_i+\frac{m^2(m+2)^2}{3(n+2)(n+4)}\right),\\
&\beta_2=\frac{m(m+2)}{3}(1-2\frac{m+2}{n+4}+\frac{m(m+2)}{(n+2)(n+4)})
\end{align*}

where $\sum y_i=\sum_{1\leq i\leq m} y_i$, $\sum y_i^2=\sum_{1\leq i\leq m} y_i^2$, 
$\sum y_iy_j=\sum_{1\leq i<j\leq m} y_iy_j$.

\medskip

\begin{theorem}[\cite{BCN}]\label{th1}
Let ${\mathcal D}\subset \Gmn$ be a finite set. Then,
\begin{enumerate}
\item for all $\mu$, $\sum_{p,p'\in {\mathcal D}}P_{\mu}(y_1(p,p'), \dots, y_m(p,p'))\geq 0$.
\item The set ${\mathcal D}\subset \Gmn$ is a $2t$-design if and only if 
for all $\mu$, \\
$1\leq \deg(\mu)\leq t$, 
$\sum_{p,p'\in {\mathcal D}}P_{\mu}(y_1(p,p'), \dots, y_m(p,p'))=0$.
\end{enumerate}
\end{theorem}

\begin{remark} The first property is basic to the so-called
{\em linear programming method} to derive bounds for codes and designs
(see \cite{B}).
\end{remark}

\subsection{Some subsets of $\Gmn$ associated with a lattice.}

Let $L\subset \R^n$ be a lattice. We define certain natural finite subsets
of $\Gmn$  
associated with $L$, in the following way. Let $\Sm$, $\Smd$, $\Smsd$
be the spaces of
$m\times m$ real symmetric, respectively real  positive definite,
and real positive semi-definite matrices.

\begin{definition}
Let $S\in \Smd$.
Let $L_S$ be the set of $p\in \Gmn$ such that $p\cap L$ is a lattice,
having a basis $(v_1,\dots, v_m)$ with $v_i\cdot v_j = S_{i,j}$ for all
$1\leq i,j\leq m$.
\end{definition}

Clearly, the sets $L_S$ are finite sets.
In the case $m=1$, the sets $L_S$ are the sets of lines supporting the 
 lattice 
vectors of fixed norm. 

\begin{definition}
Let $\delta_m(L):=\min_{S\in \Smd \mid L_S\neq \emptyset}\det S$. 
Let $S_m(L):=\cup L_S$, where 
$S\in \Smd$ and $\det S=\delta_m(L)$. The finite set $S_m(L)$ is called
the set of minimal $m$-sections
of the lattice $L$.
\end{definition}
In
par\-ti\-cu\-lar, $\delta_1(L)=\min(L)$. 
The minimal $1$-sections are the lines supporting the minimal
vectors of the lattice.

\section{Grassmannian designs and Rankin constants of lattices}

Beside the classical Hermite function $\gamma$ ($=\gamma_1$ in what follows), 
Rankin defined a collection of functions $\gamma_{m}$ associated with a lattice 
$L\subset \R^n$:

\begin{equation}
\gamma_m(L) := \delta_{m}(L)/(\det L)^{\frac{m}{n}}
\end{equation}
Thus, for $m=1$, $\gamma_1(L)$ is the classical Hermite invariant of
$L$. As a function on the set of $n$-dimensional positive definite
lattices, $\gamma_m$ is
bounded, and the supremum, which actually is a maximum, is denoted by
$\gamma_{m,n}$.  In \cite{C}, a 
characterization of the local maxima of $\gamma_m$ 
 was given.
 
\begin{definition}
\begin{enumerate}
\item A lattice $L$ is called $m$-perfect if the endomorphisms $\pr_p$
when $p\in S_m(L)$ generate $\End^s(E)$
\item A lattice $L$ is $m$-eutactic if there exist positive
coefficients $\lambda_p$, $p \in S_m(L)$ such that $\sum_{p \in \mathcal S_m (L)}
\lambda_p \pr_p=Id$.
\item A lattice $L$ is called $m$-extreme, if 
$\gamma_m$ achieves a local maximum at $L$.
\end{enumerate}
\end{definition}

\begin{theorem}[\cite{C}]
$L$ is $m$-extreme if and only if $L$ is both $m$-perfect and $m$-eutactic.
\end{theorem}

\begin{theorem}[\cite{V1}, \cite{BCN}]
If $S_m(L)$ is a $4$-design in $\Gmn$,  then it is $m$-extreme,
i.e. it achieves a local maximum of the Rankin function $\gamma_m$.
\end{theorem}

Following B. Venkov, who calls {\em strongly perfect} a lattice 
for which $S(L)$ is a $4$-design, we call $m$-strongly perfect
a lattice $L$ for which $S_m(L)$ is a $4$-design in $\Gmn$.
It is worth noticing that, since the number of classes of $m$-perfect
lattices is finite, the number of classes of strongly $m$-perfect
lattices is also finite.

\medskip

\noindent{\bf Examples:} The main sources of examples are the following:

\begin{itemize}
\item Small dimensional lattices gave the first examples of $m$-strongly
lattices: in that case, it can be checked directly,
using Theorem \ref{th1}. It was natural to look among the 
strongly perfect lattices, which  have been 
classified up to dimension $n\leq 12$(\cite{V1}, \cite{NV}, \cite{V2}).
These are: $A_2$, $D_4$, $E_6$, $E_7$, $E_8$, $K'_{10}$, ${K'_{10}}^*$,
 $K_{12}$. They are 
$m$-strongly perfect for all $m$, except $K'_{10}$, its dual, and $K_{12}$,
which are only 
$1$-strongly perfect. 

\item Extremal modular lattices. In that case, the spherical theta series
of the lattices can be used to prove strong perfection. This argument
generalizes in principle to $m>1$. Only for $m=2$ and the even unimodular
case explicit calculations on the spaces of vector-valued Siegel modular forms
show that certain families of lattices are $2$-strongly perfect,
namely the extremal ones of dimension $32$ and $48$ (see \cite{V1}, \cite{BV},
\cite{BN}).

\item Lattices with an automorphism group whose natural representation
satisfies the criterion of Theorem \ref{orb} of the next section. 
This case leads to many examples (see Table 1),
and to an infinite family of $m$-strongly perfect
lattices: the sequence of the Barnes-Wall lattices, which will be discussed in section 5.
\end{itemize}

\section{Orbits of finite subgroups of $\On$.}

A natural way to produce finite subsets of $\Gmn$ is to take the orbit
of a point under the action of a finite subgroup $G$ of $\On$.
In \cite{BCN}, we prove a criterion on the representations of $G$ 
for these sets to be designs, which naturally extends a well-known criterion
for the spherical designs.

\begin{theorem}[\cite{BCN}]\label{orb} Let $m_0\leq n/2$. Let $G<\On$ be a 
finite group. The following conditions are 
equivalent:

\begin{itemize}
\item For all $m\leq m_0$ and all $p\in \Gmn$, $G.p$ is a $2t$-design
\item For all $\mu$, $1\leq \deg(\mu)\leq t$, $l(\mu)\leq m_0$,
$(V_n^{2\mu})^G=\{0\}$
\end{itemize}
\end{theorem}

\proof We give here a simplified proof. 
Assume $\D=G.p$ is the orbit of $p\in \Gmn$. Let $G_p$ be
the stabilizer of $p$. Then,

\begin{align*}
\sum_{x\in \D}f(x)&=\frac{1}{|G_p|}\sum_{g\in G}f(g.p)\\
&=\frac{1}{|G_p|}\sum_{g\in G}(g^{-1}.f)(p)\\
&=\frac{|G|}{|G_p|}(\epsilon_G.f)(p)
\end{align*}
where $\epsilon_G=\frac{1}{|G|}\sum_{g\in G}g$.
The condition $(V_n^{2\mu})^G=\{0\}$ is equivalent to $\epsilon_G(V_n^{2\mu})=\{0\}$ which from previous equalities and the characteristic condition
(\ref{e2}) lead to the statement.

\eproof

\noindent{\bf Examples:} It is well-known that
the Weyl groups of irreducible root systems $W(R)$ 
acting on the space of homogeneous polynomials of degree $2$ 
leave invariant only the quadratic form $x_1^2+x_2^2+\dots+x_n^2$.
Therefore, these groups give rise to $2$-designs 
on all the Grassmannian spaces.
Moreover, the property for the degree $4$ holds also for 
$A_2$, $D_4$, $E_6$, $E_7$ and the degree $6$ is fulfilled for $E_8$.
It is easily checked directly on the groups; note that the partitions to be
taken into account are not only $(4)$ and $(6)$ but also, when $n\geq 4$
$(2,2)$, $(4,2)$, and, when $n\geq 6$, $(2,2,2)$.

The group $2.Co_1$ has the required property for the degree $10$, 
with no restriction on $m$.

Another interesting example is the sequence of real Clifford groups 
${\mathcal C}_k$ which are subgroups of $O(\R^{2^k})$, 
leading to $6$-designs in all the Grassmannians.
Next section considers this group and 
one subgroup of index $2$ which is the automorphism 
group of the Barnes-Wall lattice.

When the previous theorem can be applied to the group of automorphisms
of a lattice $L$, since obviously the sets $L_S$ are unions of orbits 
under the action of $\Aut(L)$, we obtain that all these sets are designs.

When the strength is equal to $4$, the possible partitions are
$(2)$, $(4)$, $(2,2)$. We have investigated the behavior of 
$\Aut(L)$ for all the lattices $L$ of dimension $4\leq n\leq 26$ which are 
known to be strongly perfect. The results are summarized in Table 1,
where only
one lattice among $\{L,L^*\}$ appears, even when they are not similar lattices.

The following situations occur (encoded in the last column of 
the table):

\begin{enumerate}
\item $G=\Aut(L)$ satisfies $(V_n^{\mu})^G=\{0\}$ for the three possible 
partitions $(2)$, $(4)$, $(2,2)$. In that case, the sets $L_S$ are $4$-designs
for all $S$, and in particular $L$ is strongly $m$-perfect for all $m$.
It holds also for any lattice with the same automorphism group, 
especially for the dual lattice.

\item $G=\Aut(L)$ satisfies $(V_n^{\mu})^G=\{0\}$ only for
 $(2)$ and  $(4)$. We can only conclude that 
the sets $L_m:=\{x\in L\mid x\cdot x=m\}$, also called the {\em layers}
of the lattice are $4$-designs, as well as the layers of the dual lattice.

\item $G=\Aut(L)$ does not satisfy $(V_n^{\mu})^G=\{0\}$ for $(2)$ and  $(4)$.

\end{enumerate}

Moreover, one can ask if any of these lattices have an automorphism group
holding the property of Theorem \ref{orb} for $t\geq 3$. It is well-known for
the Leech lattice and $t=5$ (and not for $t=6$), and next section proves 
that the lattices $E_8$ and  $\Lambda_{16}$ reache $t=3$. A direct calculation
shows that the minimal vectors of $E_8$ and $\Lambda_{16}$ do not hold an $8$-design, so
$t=3$ is the maximum. The classification of
the  integral lattices of minimum $m\leq 5$ whose set of minimal vectors
is a $6$-design, performed in \cite{Ma}, 
shows that the other lattices in this table cannot exceed $t=2$, except
possibly the lattice $N_{16}$, and the lattices $O_{23}$ and $\Lambda_{23}$
(the lattice $O_{23}$ is missing in the list of lattices given in 
\cite[Th\'eor\`eme]{Ma}, see the Erratum at: http://www.math.u-bordeaux1.fr/~martinet).
A direct computation on the automorphism groups
shows that $t=3$ is the maximum for $O_{23}$ and $\Lambda_{23}$, respectively
 $t=2$ for $N_{16}$.

The list of these lattices is taken from
\cite{V1}, with an additionnal lattice of dimension $26$ which was 
pointed to me by J. Martinet (named $T26$ after \cite{NS}. The lattice
$N_{26}$ appears in \cite{NS} as $Beis_{26}$ and $S_6(3)C_3.2$.)

We have kepted the notations of \cite{V1} for the names of the lattices,
except of course for the last one. 
The determinant is given 
in the third column, in a form that reveals the structure of
the discriminant group $L^*/L$. The automorphism group 
is given in the fifth column, with the notations of
\cite{NP}, \cite{NS}  when available.
In \cite{V1} and \cite{NS} 
more informations on these lattices are given.

The condition on $(V_n^{\mu})^G$ is checked using the Schur polynomials
associated with $\mu$.

A completely different reason for the existence of spherical designs 
in  lattices is often given by  the theory of modular forms
(see \cite{V1}, \cite{BV}). Among the list of Table 1, only
the $21$-dimensional lattice escapes from both the group theory argument
and the modular forms argument.
It is worth pointing out that it is the only one of which the dual lattice 
does not have a $4$-spherical design on its minimal vectors. Of course,
it is expected that the situation is completely different when the 
dimension grows, and the above list is anyway complete only up to
the dimension $12$.

\begin{table}\label{table}
\caption{}
$$
\begin{array}{|c|c|c|c|c|c|}
\hline
\text{dim}& \text{name}& \text{det}& \text{min} &\text{G}&\text{case}\\
\hline
4&D_4&4&2&W(F_4) & (1)\\
\hline
6 & E_6 & 3&2&2\times W(E_6) &(1)\\
\hline
7&E_7&2&2&W(E_7)&(1)\\
\hline
8 & E_8 & 1 & 2&W(E_8) & (1)\\
\hline
10& K'_{10}&6^2\cdot 3^3&4&(6\times SU(4,2)):2&(2)\\
\hline
12& K_{12}& 3^6&4&6.SU_4(3).2^2 &(2)\\
\hline
14&Q_{14}& 3^7& 4&2\times G_2(3)&(1)\\
\hline
16 & \Lambda_{16} & 2^8&4& 2^9_+\Omega^+(8,2)&(1)\\
\hline
-&O_{16}&2^6&3&D_8^4.S_6(2) &(1)\\
\hline
-&N_{16}&5^8&6&2.Alt_{10}&(2)\\
\hline
18&K'_{18}&3^5&4& (2\times 3^{1+4}:Sp_4(3)).2&(2)\\
\hline
20&N_{20}&2^{10}&4&(SU_5(2)\times SL_2(3)).2&(2)\\
\hline
-&N_{20}'&-&-&2.M_{12}.2&(2)\\
\hline
-&N_{20}''&-&-&HS_{20}&(3)\\
\hline
21& K'_{21}&12\cdot 3&4&2^{11}.3^6.5.7&(3)\\
\hline
22& O_{22}&3&3& [\Aut(\Lambda_{22}):\Aut(O_{22})]=3&(1)\\
\hline
-& \Lambda_{22}&6\cdot2 &4& (2\times PSU_6(2)).S_3&(1)\\
-&\Lambda_{22}[2]&6\cdot 2^{19}&6&-&-\\
\hline
-& M_{22}&15 &4& (2\times McL).2&(1)\\
-&M_{22}[5]&15\cdot 3^{20}&10&-&-\\
\hline
23&O_{23}&1&3&2\times CO_2&(1)\\
- &\Lambda_{23}&4&4&-&-\\
\hline
-&M_{23}&6&4&2\times CO_3&(1)\\
-&M_{23}[2]&6\cdot3^{21}&10&-&-\\
\hline
24&\Lambda_{24}&1&4&2.CO_1&(1)\\
\hline
24&N_{24}&3^{12}&6&SL_2(13)\circ SL_2(3)&(3)\\
\hline
26&N_{26}&3^{13}&6&S_6(3)C_3.2&(3)\\
\hline
26&T_{26}&3&4&3D_4(2):3&(3)\\
\hline
\end{array}
$$
\end{table}

\newpage

\section{The group $\Aut(BW_n)$}

In this section we study the tensor invariants of the automorphism group of
the Barnes-Wall lattices. We shall make use of the methods and results 
developed in \cite{NRS}. Let us recall from \cite{NRS}
some facts about the Clifford groups $\Ck$ and the Barnes-Wall lattices. 

We set $n=2^k$. The real space $\R^n$ is endowed with an orthonormal
basis $(e_u)_{u\in F_2^k}$ indexed by the elements of $\F_2^k$.

The Barnes-Wall lattice $BW_n\subset \R^n$ is the lattice defined by:

\begin{equation*}
BW_n=<2^{\lfloor \frac{k-d+1}{2}\rfloor}\sum_{u\in U}e_u, U>_{\Z}
\end{equation*}
where $U$ runs over all affine subspaces of $F_2^k$, and $d=\dim(U)$.

The first lattices of the sequence are well-known:
$BW_4\simeq D_4$, $BW_8\simeq E_8$, $BW_{16}\simeq \Lambda_{16}$
the laminated lattice of the dimension $16$.
Suitably rescaled, $\min(BW_n)=2^{\lfloor \frac{k}{2}\rfloor}$, and
$BW_n$ is even unimodular when $k\equiv 1 \mod 2$, respectively $2$-modular 
when $k\equiv 0 \mod 2$.

Bolt, Room and Wall (\cite{BRW1}, \cite{BRW2}, \cite{Bo}) and later
Brou\'e-Enguehard \cite{BE} described $\Aut(BW_n)$. 
When $n\neq 8$, it is a subgroup of index $2$ in 
the Clifford group $\Ck$ which we describe now.

The extra-special $2$-group $2^{1+2k}_+$ has a representation $E$
in $\R^n$: if 

\begin{equation*}
X(a): e_u\to e_{u+a}\text{ and }
Y(b):e_u\to (-1)^{b\cdot u}e_u,
\end{equation*}

\begin{equation*}
E=<-I,X(a),Y(b) \mid a,b\in \F_2^k>.
\end{equation*}

\begin{definition} The Clifford group ${\mathcal C}_k$ is the normalizer 
in $O(\R^n)$  of $E$.
\end{definition}

Since $q(x):=x^2$ defines a quadratic form on $E/Z(E)\simeq \F_2^{2k}$,
non degenerate and of maximal Witt index, and since
${\mathcal C}_k$ acts on $E$ (by conjugation) preserving $q$,
it induces a subgroup of $O^+(2k,2)$. 
It turns out that
the whole $O^+(2k,2)$ is realized, yielding the isomorphism:

\begin{equation*}
\Ck\simeq 2^{1+2k}_+. O^+(2k,2)
\end{equation*}

The group $O^+(2k,2)$ has a unique subgroup of index $2$, $\Omega^+(2k,2)$.
Its parabolic subgroups are the stabilizers in  $\Omega^+(2k,2)$ of totally
isotropic subspaces; they are maximal in $\Omega^+(2k,2)$.
Let $P(2k,2)$ be the one associated with the image in $\F_2^{2k}$ 
of $<\pm X(a) \mid a \in \F_2^k>$.

According to  \cite{NRS}, the following transformations 
are explicit generators of
the group $ {\mathcal C}_k$:

\begin{enumerate}

\item Diagonal transformations: $ e_u\to (-1)^{q(u)}e_u$, where $q$ is any
binary quadratic form, and $-I$.

\item Permutation transformations: $e_u\to e_{\phi(u)}$, where $\phi\in AGL(k,2)$.

\item $ H:=h\otimes I_2\otimes\dots\otimes I_2$, 
$ h=\frac{1}{\sqrt{2}}\begin{bmatrix}1&1\\1&-1\end{bmatrix}$
(here $\R^n$ and $(\R^2)^{\otimes k}$ are identified in an obvious way).
\end{enumerate}

Straightforward calculations show that these elements normalize $E$.
Moreover, the induced action of the elements of the first and second type
on $\F_2^{2k}$ is given by the respective matrices
 $\begin{bmatrix} 1&b\\0&1\end{bmatrix}$ where $b$ is
the symplectic matrix associated with $q$, and 
$\begin{bmatrix} \phi &0\\0&\phi^{-tr}\end{bmatrix}$
where $\phi\in GL(2,k)$. The group generated by these transformations
on $\F_2^{2k}$ is the parabolic group $P(2k,2)$.

The element $H_2:=h\otimes h\otimes  I_2\otimes\dots\otimes I_2$
has rational entries. The subgroup $\Gk$ of $\Ck$ generated
by the elements of the first and second type, and $H_2$,
generate a subgroup of $\Omega^+(2k,2)$, containing $P(2k,2)$,
hence equal to $\Omega^+(2k,2)$. It follows that $\Gk$ has index $2$ in $\Ck$
and is rational; hence it is the automorphism group of $BW_n$ (see \cite{NRS}). 

\smallskip

The polynomial invariants of $\Ck$ are described, first by B. Runge
(\cite{R1}, \cite{R2}, \cite{R3}), then with a different proof by G. Nebe, E. Rains, N.J.A. Sloane 
(\cite{NRS}, in terms of self-dual binary codes. As a consequence,
the first non trivial invariant occurs for the degree $8$, associated with 
the first non trivial self-dual binary code which is the $[8,4,4]$ 
Hamming code. We extend here this result to the subgroup $\Gk$.

\begin{theorem}\label{th} If $k\geq 3$ and $d\leq 6$, then 

$$(V^{\otimes d})^{\Gk}=(V^{\otimes d})^{\Ck}=(V^{\otimes d})^{O(V)}.$$
\end{theorem}

\begin{corollary} The orbits of $\Aut(BW_n)$ on $\Gmn$ are $6$-designs.
In particular, the sets $(BW_n)_S$ are $6$-designs and the lattice $BW_n$
is strongly $m$-perfect for all $m$.
\end{corollary}

\begin{remark}

- Theorem \ref{th} shows more than what is needed for the Grassmannian design
property,  since $V^{\otimes 6}$ contains the representations
associated with arbitrary partitions of degree lower or equal to $6$.

- The fact that the set of minimal vectors is a $6$-spherical design 
was already proved by direct calculation by Boris Venkov (\cite{V1}).

\end{remark}

\proof  The argument in \cite{NRS} extends straightforwardly to
the {\em tensor invariants} of $\Ck$. 
Let $V:=\R^n$. To a binary code $C$ of length $d$,
is associated a {\em tensor enumerator} $T_C^{(k)}\in V^{\otimes d}$.
To a $k$-tuple $(w_1,\dots,w_k)$ of codewords, we associate a 
$k\times d$
matrix which rows are the words $w_1,\dots,w_k$. Let 
$u_1,\dots,u_d$ be the $d$ columns of this matrix. Then:

\begin{equation*}
T_C^{(k)}:=\sum_{(w_1,\dots,w_k)\in C^k} e_{u_1}\otimes\dots \otimes e_{u_d}
\end{equation*}
where 

The  usual (generalized) weight enumerator $W_C^{(k)}$ is
obtained by the symmetrization $V^{\otimes d}\to {\operatorname{Sym}}^d(V)$.
For the same reasons, when $C$ is self-dual, $T_C^{(k)}$ is invariant
under the action of $\Ck$. A straightforward generalization of the proof
in \cite{NRS} of the fact that the invariants of $\Ck$ 
on ${\operatorname{Sym}}^d(V)$ are exactly spanned by the polynomials 
$W_C^{(k)}$ when $C=C^{\perp}$ shows that the invariants of $\Ck$ 
on $V^{\otimes d}$ are spanned by the $T_C^{(k)}$ when $C=C^\perp$.
To determine the invariants of $\Gk$, we follow the same steps 
as in \cite{NRS}: the first is the description of 
$(V^{\otimes d})^{\Pk}$, which we recall in next lemma.

\begin{lemma}[\cite{NRS}, Theorem 4.6]\label{l1}
The space $(V^{\otimes d})^{\Pk}$ is generated by the $T_C^{(k)}$
where $C$ runs over the binary codes of length $d$ 
such that ${\bf 1}\subset C\subset C^{\perp}$ and $\dim(C)\leq k+1$.
\end{lemma}

The second calculates $\epsilon_{\Pk}H_2$
as a linear combination of the $T_C^{(k)}$ associated with binary codes 
satisfying ${\bf 1}\subset C\subset C^{\perp}$ (which obviously belong to 
$(V^{\otimes d})^{\Pk}$; only those with 
$\dim(C)\leq k+1$ are linearly independent).

\begin{lemma}\label{l2} Let $C$ be a  binary code of length $d$ 
such that ${\bf 1}\subset C\subset C^{\perp}$ and $\dim(C)\leq k+1$. 
Let $r:=d/2-\dim(C)$. 

$$(\epsilon_{\Pk}H_2). T_C^{(k)}=a_1T_C^{(k)}+
a_2\sum_{\substack{C'\subset {C'}^\perp\\ C\subset C', [C':C]=2}} T_{C'}^{(k)} +
a_4\sum_{\substack{C'\subset {C'}^\perp\\ C\subset C', [C':C]=4}}  T_{C'}^{(k)}
$$
where
$$
\begin{cases}
a_1=2^{-2r}(1+2\frac{(2^{2r}-1)(2^{2r-2}-1)}{(2^{k}-1)(2^{k-1}-1)}
-3\frac{2^{2r}-1}{2^k-1})\\
a_2=\frac{3.2^{-2r}}{2^k-1}(1-\frac{2^{2r-2}-1}{2^{k-1}-1})\\
a_4=\frac{3.2^{-2r}}{(2^k-1)(2^{k-1}-1)}
\end{cases}
$$

Moreover, $a_1=1$ if and only if $r=0$ or $r=k$.
\end{lemma}

\proof Let $\mu(w_1,\dots,w_k):=e_{u_1}\otimes\dots\otimes e_{u_d}$.
We have (as a consequence of the Poisson summation formula)

\begin{equation*}
H_2T_C^{(k)}=2^{-2r}\sum_{\substack{w_1,w_2\in C^{\perp}\\
w_3,\dots,w_k\in C}} \mu(w_1,\dots,w_k)
\end{equation*}

As a consequence of the change from $H$ to $H_2$, not only the first, but also
the second vector is allowed to be in $C^\perp$. Therefore, by the same
argument as in \cite{NRS}, there exists coefficients $a_1,a_2,a_4$ (depending on $r$ and $k$)
 such that 

\begin{equation}\label{e3}
\epsilon_{\Pk}H_2 T_C^{(k)}=a_1T_C^{(k)}+
a_2\sum_{\substack{C'\subset {C'}^\perp\\ C\subset C', [C':C]=2}} T_{C'}^{(k)} +
a_4\sum_{\substack{C'\subset {C'}^\perp\\ C\subset C', [C':C]=4}}  T_{C'}^{(k)}
\end{equation}
and we are left with the computation of these coefficients. 
Let $<,>$ denote the scalar product induced on $V^{\otimes d}$
by the  Euclidean structure on $V$. For any 
codes $C$, $D$, with ${\bf 1}\subset C \subset D \subset D^\perp \subset C^\perp$, we have:

\begin{equation*}
<T_C^{(k)}, T_D^{(k)}>=|C|^k
\end{equation*}
and 

\begin{equation*}
<(\epsilon_{\Pk}H_2).T_C^{(k)}, T_D^{(k)}>=<H_2.T_C^{(k)}, T_D^{(k)}>=2^{-2r}[D:C]^2|C|^k.
\end{equation*}

Let $n_2^r$, respectively $n_4^r$ be
 the number of self-orthogonal codes containing 
$C$ to index $2$, respectively  $4$. Obviously, $n_2^r$ equals the number
 of isotropic lines in the symplectic space $C^\perp/C$ of dimension $2r$,
and $n_4^r$ equals the number
 of totally isotropic planes  in $C^\perp/C$. Therefore,
$n_2^r=2^{2r}-1$ and $n_4^r=(2^{2r}-1)(2^{2r-2}-1)/3$.
Taking the scalar product of equation (\ref{e3}) with $T_D^{(k)}$,
successively for $D=C$, then for a self-orthogonal code containing $C$
to index $2$ and $4$, we obtain the three equations (after having divided by 
$|C|^k$):

\begin{align*}
2^{-2r}&=a_1+a_2n_2^r+a_4n_4^r\\
2^{-2r}.4&=a_1+a_2.2^k+a_2(n_2^r-1)+a_4n_2^{r-1}.2^k+a_4(n_4^r-n_2^{r-1})\\
2^{-2r}.16&=a_1+a_2.3.2^k+a_2(n_2^r-3)\\
          &\quad +a_4.4^k+a_4(3n_2^{r-1}-3).2^k+a_4(n_4^r-3n_2^{r-1}+2)
\end{align*}

which lead to the expressions of Theorem \ref{th}.

\eproof

We end the proof of the theorem in the same way as in \cite{NRS}.
We have $(V^{\otimes d})^{\Gk}=\ker(\epsilon_{\Pk}H_2-I)\cap (V^{\otimes d})^{\Pk}$. From Lemma \ref{l2}, when the elements $T_C^{(k)}$
are ordered by increasing $\dim(C)$, 
the matrix of the  transformation $\epsilon_{\Pk}H_2$ is upper triangular.
If $k\leq 3$ and $d\leq 6$, the only diagonal coefficients which are equal to
$1$ correspond to $C=C^\perp$ and we can conclude by \cite{NRS}, Lemma 4.8

\eproof

\begin{remark} Of course, for arbitrary degree $d$, the group $\Gk$ has 
more invariants than $\Ck$. For $k=2$ and $d=6$, 
we have $a_1=1$ for $r=2$, i.e. for the code $C=\bf{1}$. The element

$$T_{\bf{1}}^{(2)}-\frac{1}{12}\sum_{\substack{{\bf 1}\subset C\subset C^\perp\\\dim(C)=2}}T_{C}^{(2)}$$
is the unique degree $6$ additional  invariant under ${\mathcal G}_2$.

For $k=3$ and $d=8$, the situation is the same, with 

$$T_{\bf{1}}^{(3)}
-\frac{1}{40}
\sum_{\substack{{\bf 1}\subset C\subset C^\perp\\\dim(C)=2}}T_{C}^{(3)}
+\frac{1}{480}
\sum_{\substack{{\bf 1}\subset C\subset C^\perp\\\dim(C)=3}}T_{C}^{(3)}
$$
as an invariant of degree $8$. The Molien series confirms that the degree $8$
polynomial invariant space has dimension $3$, spanned by the two classes
of self-dual codes and this one.

\end{remark}

\section{Other Grassmannian designs}

When a group $G$ is known to fulfill the conditions of Theorem \ref{orb},
among its orbits the most interesting ones are  the ones shorter 
than the ``generic'' ones, i.e. the ones with a non trivial
isotropic group. In general, it is not easy to describe these orbits.
In the case of the Clifford group $\Ck$,
some of these orbits are described in a very explicit way in \cite{CHRSS},
in view of the construction of Grassmannian {\em codes} for the
chordal distance. We next discuss under which conditions certain smaller subsets
of these sets remain to be $6$-designs or $4$-designs. More precisely,
we prove that it depends on a similar condition of design associated
with the underlying finite geometry.

\subsection{The construction}

The alluded construction is the following. Let $S\subset \F_2^{2k}$ be a 
totally isotropic subspace of dimension $k-s$. The preimage $\tilde{S}$ of $S$ in $E$
is an abelian group, $2$-elementary.
(The identification between $\F_2^{2k}$ and $E/\{\pm 1\}$ is still the same, 
sending $X(a)Y(b)$ to $(a,b)$),
 It decomposes the space $V=\R^n$ into  
$2^{k-s}$ irreducible subspaces of dimension $2^s$, which are pairwise orthogonal.
Let ${\D}_S\subset \G_{2^s,2^k}$ be the set of these $2^{k-s}$ subspaces.

More generally, if $\Sigma$ is a set of isotropic subspaces
of the same dimension $k-s$, we set

\begin{equation}\label{ds}
{\D}_{\Sigma}:=\cup_{S\in \Sigma} {\D}_S \subset \G_{2^s,2^k}.
\end{equation}

\medskip

\noindent {\bf Example:} We can take $\Sigma$ to be the whole set of totally isotropic subspaces
of fixed dimension $k-s$. In that case, $\Sigma$ is a single orbit under 
$O^+(2k,2)$, and ${\D}_{\Sigma}$ is a single orbit under $\Ck$. Therefore it is a $6$-design.
Note that, when $s=0$, $\Sigma$ splits into two orbits under the action of
$\Omega^+(2k,2)$; the set of lines corresponding to one orbit is the
set of lines supporting the minimal vectors of $BW_n$, $n=2^k$.

\medskip
Of course, we are interested in the smallest possible sets, and the above 
example is the largest one.
A natural question is then: which conditions should $\Sigma$ satisfy, so that
$\Dsig$ is a design? How small can we take $\Sigma$?
To answer these questions, we need two more ingredients:
another criterion for Grassmannian designs, and the notion of designs on the spaces
of totally isotropic subspaces of fixed dimension.

\subsection{A new criterion} 

Let $\D\subset \Gmn$. Let $\sigma:=y_1+y_2+\dots +y_m$.

\begin{theorem}\label{crit}
For all $m,n$, and $t$, there exists a constant $c_{m,n}(2t)$ such that
\begin{enumerate}
\item For all $\D\subset \Gmn$, $\frac{1}{|\D|^2}\sum_{p,p'\in \D} \sigma(p,p')^t\geq c_{m,n}(2t)$.
\item $\D
$ is a $2t$-design if and only if  $\frac{1}{|\D|^2}\sum_{p,p'\in \D} \sigma(p,p')^t= c_{m,n}(2t)$.
\end{enumerate}
\end{theorem}

\proof From the defining property of Grassmannian designs (Definition \ref{def1}), since $\sigma^t$
has degree $t$ in the variables $y_1,\dots,y_m$, if $\D$ is a $2t$-design,

$$\frac{1}{|\D|^2}\sum_{p,p'\in \D} \sigma(p,p')^t= \int_{[0,1]^m}\sigma^t d\mu(y_1,\dots,y_m)$$

We set $c_{m,n}(2t):=\int_{[0,1]^m}\sigma^t d\mu(y_1,\dots,y_m)$.

\begin{lemma}\label{l3} There exists positive coefficients $\lambda_{t,\mu}>0$  such that:

\begin{equation*}
\sigma^t=\sum_{\mu, \deg(\mu)\leq t}\lambda_{t,\mu}P_\mu.
\end{equation*}

\end{lemma}

\proof For $t=1$, we have $\sigma=\frac{m(n-m)}{n}P_{(1)}+\frac{m^2}{n}$.

For $t>1$, we proceed by induction. Let us assume first that $\deg(\mu)<t$.
We have 

\begin{align*}
[\sigma^t,P_{\mu}]&=[\sigma^{t-1},\sigma P_{\mu}]\\
                  &=[\sigma^{t-1},(\frac{m(n-m)}{n}P_{(1)}+\frac{m^2}{n})P_{\mu}]
\end{align*}
We know  that $P_{(1)}P_{\mu}$ is a 
linear combination with non negative  coefficients of the $P_{\kappa}$ 
(\cite[Lemma 2]{BBC}). By induction, we obtain

\begin{align*}
[\sigma^t,P_{\mu}]&\geq [\sigma^{t-1},\frac{m^2}{n}P_{\mu}]\\
                  &> 0 \text{ (if } \deg(\mu)<t).
\end{align*}

If $\deg(\mu)=t$, the second  term of the first inequality is zero. We need more information on the expression $\sigma P_{\mu}$ on the $P_{\kappa}$.
The analogue of the ``three-term relation'' 
for  orthogonal polynomials in one variable gives (see \cite{B}):

\begin{equation*} 
\sigma P_{\mu}=\sum_{ \deg(\kappa)=k+1}A_k[\mu,\kappa]P_{\kappa}+
\sum_{ \deg(\kappa)=k}B_k[\mu,\kappa]P_{\kappa}+
\sum_{ \deg(\kappa)=k-1}C_k[\mu,\kappa]P_{\kappa}
\end{equation*}
where $k=\deg(\mu)$. Moreover, 
$C_k[\mu,\kappa][P_{\kappa},P_{\kappa}]=A_{k-1}[\kappa,\mu][P_{\mu},P_{\mu}]$ is zero unless
$\mu$ is obtained from $\kappa$ by the increase of one of its parts by one,
in which case $A_{k-1}[\kappa,\mu]>0$ (\cite{B}).
In $[\sigma^t,P_{\mu}]=[\sigma^{t-1},\sigma P_{\mu}]$ only those terms (and at least one)
give a contribution, so, by induction, we obtain the desired property.

\eproof

Since $c_{m,n}(2t)=[\sigma^t,1]=\lambda_{t,0}$,
and from the design criterion and the positivity condition of Theorem \ref{th1}, 
the proof of Theorem \ref{crit} is completed (note that it is crucial than none of
the $\lambda_{t,\mu}$ is equal to zero).

\eproof

\begin{remark} This criterion is analogous to \cite[Th\'eor\`eme 8.1]{V1}, and similar versions 
exist in principle for any notion of design. We shall come across a similar criterion for the 
designs of totally isotropic spaces.
\end{remark}

\medskip

\begin{remark}\label{r1}
It is not apparently  easy to calculate 
$c_{m,n}(2t)$ by the integration  formula 
$c_{m,n}(2t):=\int_{[0,1]^m}\sigma^t d\mu(y_1,\dots,y_m)$. 
It is worth noticing that, since 
$c_{m,n}(2t)=[\sigma^t,1]=\lambda_{t,0}$, it becomes easy once
one has calculated explicitly the polynomials $P_{\mu}$ for $\deg(\mu)\leq t$.
For example, we obtain from \S 2.1,

\begin{align*}
c_{m,n}(2)&= \frac{m^2}{n}\\
c_{m,n}(4)&=\frac{m^2}{3n}\left(\frac{2(m-1)^2}{n-1}+\frac{(m+2)^2}{n+2}\right)
\end{align*}
and, using $[\sigma^3,1]=[\sigma^2,\sigma]$ and 
$[P_{(1)},P_{(1)}]=\dim(V_n^{(2)})^{-1}=\frac{2}{(n-1)(n+2)}$,
we can even calculate

\begin{align*}
c_{m,n}(6)&= \frac{m^2}{3n}
\left(\frac{(m-1)^2(m+2)^2}{(n-1)(n+2)}(\frac{2n}{n-2}+\frac{n+3}{n+4})\right.\\
    &\qquad\qquad\qquad\left.-8\frac{m(m-1)^2}{(n-1)(n-2)}+
\frac{(m+2)^2(2m+3)}{(n+2)(n+4)}\right)
\end{align*}
\end{remark}

\subsection{The space of totally isotropic subspaces}

Let $X_{w}$ be the set of totally isotropic subspaces of dimension $w\leq k$
of the quadratic space $(\F_2^{2k},q)$. 
The group $G:=O^+(2k,2)$ acts transitively
on $X_w$; the stabilizer of an element is a maximal parabolic subgroup $P_w$.
The orbits of $G$ on pairs of elements $(S,S')$ (also called orbitals)
are investigated in \cite{WW}; they are characterized by two quantities:
$\dim(S\cap S')$ and $\dim(S \cap S'^\perp)$. Since 
$\dim(S \cap S'^\perp)=\dim(S^\perp \cap S')$, they are symmetric.
In the special case $w=k$ of the maximal totally isotropic subspaces,
$S=S^\perp$ and one value is enough, giving to $X_w$ the structure 
of a $2$-point homogeneous space 
(for the distance $d(S,S')=k-\dim(S \cap S')$).

The space $L(X_w)$ of complex valued functions on $X_w$ decomposes, under
the action of $G$, into irreducible subspaces with multiplicities 
equal to one; to each subspace is associated a unique zonal function.
This decomposition, and the corresponding zonal functions, are computed in
\cite{S1}, \cite{S2} (in \cite{S2}, the general case of Chevalley
groups over $\F_q$ is treated; it is assumed that the characteristic
is different from $2$, although the situation would be completely analogous.
In \cite{S2}, the case $w=k$ is treated in full generality).
We only need here the general form  of this decomposition
(\cite[Theorem 6.23]{S1}):

\begin{equation*}
L(X_w)=\oplus_{(m,r)\in I}V_{m,r}
\end{equation*}
where $I:=\{(m,r)\mid 0\leq m\leq w, 0\leq r\leq m\wedge (k-w)\}$.
If $y:=\dim(S\cap S')$, $x+y:=\dim(S\cap S'^\perp)$, the corresponding 
zonal (spherical) function $G_{m,r}$ is a polynomial in $2^x$, $2^y$
and of degree $m$ in $2^y$. Note that $2^y=|S\cap S'|$.
When $w=k$, these polynomials are polynomials in one variable, and
identified as $q$-Krawtchouk polynomials. In that case, 
the $t$-designs are defined in the usual way (see \cite{DL}).

\begin{theorem} 
There exists constants $d_{w,k}(t)$ such that:
\begin{enumerate}
\item For all $\Sigma\subset X_w$, 
$\frac{1}{|\Sigma|^2}\sum_{S,S'\in \Sigma}|S \cap S'|^t \geq d_{w,k}(t)$
\item When $w=k$, equality holds if and only if $\Sigma$ is a 
$t$-design.
\end{enumerate}
\end{theorem}

\begin{remark} When $w<k$, the interpretation of the case of equality in terms
of designs is not so clear. Since 
the irreducible spaces require a double index, the notion of $t$-designs itself
is not so clear, although the most natural one would be, like in the case 
of the non-binary Johnson scheme, to require 
orthogonality with $\oplus V_{m,r}$, where $(m,r)\in \{(m,r)\mid m\leq t\}$.
Then, one would need to look carefully at the positivity of the coefficients
of the expansion of $(2^y)^t$ on the $G_{m,r}$. In the case
$w=k$, the positivity is guaranteed, thanks to the three-terms relation, 
whose coefficients are the
intersection numbers of the association scheme (see 
\cite[II.2(2.1), III.1(1.2)]{BI}).

\end{remark}

Here, the computation of the numbers $d_{w,k}(t)$ is easy, since they
come from the constant term: $d_{w,k}(t)=[y^t,1]$. So,

\begin{align*}
d_{w,k}(t)&=\frac{1}{|X_w|^2}\sum_{S,S'\in X_w}|S \cap S'|^t\\
&=\frac{1}{|X_w|^2}\sum_{x,y} |\text{Orb}(x,y)|(2^y)^t
\end{align*}
where $\text{Orb}(x,y)$ is the orbital associated with the values $(x,y)$;
its cardinality is calculated in \cite[Theorem 5.5]{WW}.

\subsection{When is ${\mathcal D}_{\Sigma}$ a design?}

\begin{theorem}\label{tt} Let $\D_{\Sigma}$ be defined as in 
(\ref{ds}).

\begin{enumerate}
\item $\D_{\Sigma}$ is always a $2$-design.
\item For $t=2$ and $t=3$, $\D_{\Sigma}$ is a $2t$-design
if and only if $\Sigma$ satisfies  the equality:

$$\frac{1}{|\Sigma|^2}\sum_{S,S'\in \Sigma}|S \cap S'|^{t-1} = d_{w,k}(t-1).$$

\end{enumerate}
\end{theorem}

\proof We calculate $\frac{1}{|\D_{\Sigma}|^2}\sum_{p,p'\in \D_{\Sigma}} \sigma(p,p')^{t}$. By the construction,

\begin{equation*}
\frac{1}{|\D_{\Sigma}|^2}\sum_{p,p'\in \D_{\Sigma}} \sigma(p,p')^{t}=
\frac{1}{2^{2(k-s)}|\Sigma|^2}\sum_{S,S'\in \Sigma}(
\sum_{p\in \D_S,p'\in \D_{S'}} \sigma(p,p')^{t}).
\end{equation*}

Let $\dim(S\cap S'):=k-u$. The $2^{k-u}$ irreducible subspaces associated with
$\tilde{S\cap S'}$ are obtained from the $2^{k-s}$ ones associated with  $S$,
by summing together $2^{u-s}$ of them. These ones are precisely the ones 
on which the corresponding characters of $\tilde{S}$ and $\tilde{S'}$ 
coincide on $\tilde{S\cap S'}$.
According to \cite[(9)]{CHRSS}, if $p\in \D_S$ and
$p'\in\D_{S'}$ are contained in the same 
irreducible subspaces associated with
$\tilde{S\cap S'}$, 

$$\sigma(p,p')=2^k\frac{|S\cap S'|}{|S||S'|}=2^{2s-u},$$
and it holds for $(2^{u-s})^2$ pairs $(p,p')$.
Otherwise, $\sigma(p,p')=0$, except if $p=p'$ of course.

All together, we obtain

\begin{equation*}
\frac{1}{|\D_{\Sigma}|^2}\sum_{p,p'\in \D_{\Sigma}} \sigma(p,p')^{t}=
2^{(2s-k)t}\frac{1}{|\Sigma|^2}\sum_{S,S'\in \Sigma}|S\cap S'|^{t-1}.
\end{equation*}

From Theorem \ref{crit}, we obtain that $\D_{\Sigma}$ is a $2t$-design
if and only if 

\begin{equation}\label{eee}
\frac{1}{|\Sigma|^2}\sum_{S,S'\in \Sigma}|S\cap S'|^{t-1}=
2^{-(2s-k)t}c_{2^s,2^k}(2t).
\end{equation}

When $t=1$, from Remark \ref{r1}, $2^{-(2s-k)}c_{2^s,2^k}(2)=1$, and the
previous equality always holds.

When $t=2,3$, we know  that, taking $\Sigma=X_{k-s}$, we do obtain
a $2t$-design, and hence, that (\ref{eee}) is
fulfilled. We have proved two things:

\begin{itemize}
\item $2^{-(2s-k)t}c_{2^s,2^k}(2t)=d_{k-s,k}(t-1)$.
\item Assertion (2) of the theorem.
\end{itemize}

\eproof

It remains, of course, to give examples of sets $\Sigma$ with the
property (2) of Theorem \ref{tt}. One example is given by 
maximal spreads. These are standard objects of finite geometries.

\begin{definition} The set $\Sigma\subset X_w$ is called a spread if
$\Sigma$ is a set a totally isotropic subspaces,
such that the intersection of two distinct elements is reduced to $\{0\}$.
A maximal spread is a spread, such that $\cup_{S\in \Sigma}S$ is exactly
equal to the whole set of isotropic elements.
\end{definition}

The number of non zero isotropic vectors is $(2^k-1)(2^{k-1}+1)$.
Therefore, a maximal spread in $X_w$ must have $(2^k-1)(2^{k-1}+1)/(2^w-1)$
elements, and hence a necessary condition for the existence 
of a maximal spread, is that this number
is an integer.
It is well known that, when $w$ divides $k$, maximal spreads do exist.

\begin{theorem}\label{tt1}
Let $\Sigma$ be a maximal spread in $X_{k-s}$. Then,
$\D_{\Sigma}$ is a $4$-design.
\end{theorem}

\proof  Let $\Sigma$ be a spread, and let $N:=|\Sigma|$. We calculate 

\begin{equation*}
\frac{1}{|\Sigma|^2}\sum_{S,S'\in \Sigma}|S \cap S'| = 
\frac{1}{N^2}(N(N-1)+N.2^{k-s})=1+\frac{2^{k-s}-1}{N}.
\end{equation*}

On the other hand, from Remark \ref{r1},

\begin{equation*}
2^{-2(2s-k)}c_{2^s,2^k}(4)=\frac{2^{-(2s-k)+1}}{3}(\frac{(2^s-1)^2}{2^k-1}+
\frac{(2^{s-1}+1)^2}{2^{k-1}+1}).
\end{equation*}

The condition (\ref{eee}) leads to $N=(2^k-1)(2^{k-1}+1)/(2^{k-s}-1)$.

\eproof

\smallskip

\paragraph{\bf Aknowledgements:} I am indebted to
Eiichi Bannai, Jacques Martinet,  Gabriele Nebe, Neil Sloane
for helpful discussions and comments.

\end{document}